\documentclass[10pt,leqno,a4paper]{article}
	
	\usepackage[cp850]{inputenc}
   	\newcommand{\proof}{\noindent{\bf{Proof:}\hskip 0.5em}}
	\newcommand{\qed}{\,\,{\bf QED}}

	\newcommand{µ}{\ae}
	\newcommand{°}{\o}
	\newcommand{Õ}{\aa}
	\newcommand{ã}{\AE}
	\newcommand{Ï}{\O}
	\newcommand{┼}{\AA}
	\newcommand{÷}{\"o}
	\newcommand{Þ}{\`e}
	\newcommand{ß}{\'a}

      \newtheorem{theorem}{Theorem}[section]
      \newtheorem{proposition}[theorem]{Proposition}
      \newtheorem{corollary}[theorem]{Corollary}
      
      \newtheorem{definition}[theorem]{Definition}
      
      \newtheorem{remark}[theorem]{Remark}

\begin{document}

      \title{Operator monotone functions of several variables}
	\date{}
      \author{Frank Hansen}
      
      \maketitle

     \begin{abstract}	
     We propose a notion of operator monotonicity for functions of several variables, which
	extends the well known notion of operator monotonicity for functions of only one variable. The notion
	is chosen such that a fundamental relationship between operator convexity and operator monotonicity for
	functions of one variable is extended also to functions of several variables.
	\end{abstract}

	\section{Introduction and main result}

	The notion of operator convexity for functions of several variables has been 
	extensively studied in the literature. The first step is to define the functional calculus for functions
	of several variables. This can be done in the following way:

	Let $ I_1,\dots,I_k $ be real intervals and let
      $ f:I_1\times\cdots\times I_k\to{\mathbf R} $ be a Borel measurable and
	essentially bounded function. Let $ x=(x_1,\dots,x_k) $ be a $ k $-tuple
      of bounded self-adjoint operators on Hilbert spaces 
      $ H_1,\dots,H_k $ such that the spectrum of $ x_i $ is contained
      in $ I_i $ for $ i=1,\dots,k. $ We say that such a $ k $-tuple is in the
      domain of $ f. $ If
      $$
      x_i=\int_{I_i}\lambda_i\, E_i(d\lambda_i)\qquad i=1,\dots,k
      $$
      is the spectral decomposition of $ x_i, $ we define
      \begin{equation}\label{functional calculus}
      f(x)=\int_{I_1\times\cdots\times I_k}
      f(\lambda_1,\dots,\lambda_k)\,
      E_1(d\lambda_1)\otimes\cdots\otimes E_k(d\lambda_k)
      \end{equation}
      as a bounded self-adjoint operator on $ H_1\otimes\cdots\otimes H_k, $ 
	cf.~\cite{kn:hansen:1997:2, kn:araki:2000, kn:koranyi:1961}. If the Hilbert spaces 
	are of finite dimension, then 
	the above integrals become finite sums, and we may consider the functional calculus for 
	arbitrary real functions. This construction extends the definition of 
	Korßnyi~\cite{kn:koranyi:1961} for functions of two variables and have the 	property that 
	$$
	f(x_1,\dots,x_k)=f_1(x_1)\otimes\cdots\otimes f_k(x_k), 
	$$
	whenever $ f $ can
	be separated as a product $ f(t_1,\dots,t_k)=f_1(t_1)\cdots f_k(t_k) $ of $ k $ 
	functions each depending on only one variable.

	\begin{remark}\rm
	One might consider the functional calculus only for commuting operators 
	$ x_1,\dots,x_k $ on a single Hilbert space $ H $ and define 
	$$
	f_{com}(x_1,\dots,x_k)=\int f(\lambda_1,\dots,\lambda_k)\, dE(\lambda_1,\dots,\lambda_k)
	$$
	as an operator on $ H, $
	where $ E $ is the product measure of the commuting spectral measures associated 
	with each of the operators. This approach
	was suggested by Pedersen and Lieb in~\cite{kn:pedersen:2001,kn:Lieb:2001}. Our 
	definition in 
	equation~(\ref{functional calculus}) can then be written as
	$$
	f(x_1,\dots,x_k)=
	f_{com}(x_1\otimes 1\otimes\cdots\otimes 1,\dots,1\otimes\cdots\otimes 1\otimes x_k)
	$$
	for arbitrary non-commuting operators $ x_1,\dots,x_k $ on $ H. $ If however the 
	operators $ x_1,\dots,x_k $ do
	commute, then there is a self-adjoint projection $ P $ on $ H\otimes\cdots\otimes H $ 
	with range isomorphic to
	$ H $ such that $ f_{com}(x_1,\dots,x_k)=Pf(x_1,\dots,x_k)P. $ The two approaches are 
	thus essentially equivalent.
	\end{remark}

	Once the functional calculus is defined, we say that a function 
	$ f:I_1\times\cdots\times I_k\to{\bf R} $
	is operator convex, if $ f $ is continuous and the operator inequality
	$$
	f(\lambda x+(1-\lambda)y)\le\lambda f(x)+(1-\lambda)f(y)\qquad\forall\lambda\in[0,1]
	$$
        holds for all $ k $-tuples of self-adjoint operators
        $ x=(x_1,\dots,x_k) $ and $ y=(y_1,\dots,y_k) $ in the domain of $ f $ acting on
	 any Hilbert spaces $ H_1,\dots, H_k. $ The definition is meaningful since also the
	$ k $-tuple $ \lambda x+(1-\lambda)y $ is in the domain of $ f. $ We say that $ f $ is  
       	matrix convex of order $ (n_1,\dots,n_k), $ if the operator inequality
      	holds for operators on Hilbert spaces of finite dimensions $ (n_1,\dots,n_k). $

	The aim of this paper is to define the notion of an operator monotone function
	also for functions of several variables. The definition should, 
	when restricted to functions of only one variable, be a simple 
	reformulation of the ordinary condition for such functions.
	We also want the following theorem to be valid.

	\begin{theorem}\label{theorem:operator monotonicity and operator convexity}
	Let $ f:[0,\alpha_1[\times\cdots\times [0,\alpha_k[\to{\mathbf R} $ be a continuous real
 	function. The following statements are equivalent:

	\begin{itemize}

	\item[$(i)$] $ f $ is operator convex, and $f(r_1,\dots,r_k)\le 0 $ if $ r_i=0 $
	for some $ i=1,\dots,k. $

	\item[$(ii)$] The function
	$ g:\,]0,\alpha_1[\times\cdots\times ]0,\alpha_k[\to{\mathbf R} $ defined by setting
	$$
	g(r_1,\dots,r_k)=r_1^{-1}\cdots r_k^{-1}f(r_1,\dots,r_k)
	$$
	is operator monotone.
	
	\end{itemize}
	\end{theorem}

	The theorem above is known to be valid for functions of one 
	variable~\cite[2.4 Theorem]{kn:hansen:1982}, and
	the extension to functions of several variables seems to be very natural. 
	Our notion of operator monotonicity
	for functions of several variables is ultimately given in 
	Definition~\ref{ex post definition of operator monotonicity}, but it depends on 
	intermediary notions and results	
	given in Definition~\ref{decomposition}, Definition~\ref{definition of index}, 
	Definition~\ref{operator monotone function}, and 
	Corollary~\ref{independency of a particular index}.

	Before proceeding with this programme, we shall briefly discuss other possible 
	definitions of operator monotonicity for functions of several variables, which we 
	ultimately have rejected. 
	
	\begin{proposition} Let $ f $ be a non-negative continuous function of $ k $
	variables defined in the first quadrant $ [0,\infty[\times\cdots\times[0,\infty[. $
	If $ f $ is matrix concave of order $ (n_1,\dots,n_k), $ then 
	$$
	0\le x_i\le y_i\quad i=1,\dots,k\qquad\Rightarrow\qquad f(x)\le f(y)
	$$
	for arbitrary $ k $-tuples of positive semi-definite matrices $ x=(x_1,\dots,x_k) $ and 
	$ y=(y_1,\dots,y_k) $ of order $ (n_1,\dots,n_k). $
	\end{proposition}

	\proof Let the appropriate $ k $-tuples of matrices be chosen and take 
	$ \lambda\in[0,1[. $
	We set $ z_i=\lambda(1-\lambda)^{-1}(y_i-x_i) $ and notice that
	$$
	\lambda y_i=\lambda x_i+(1-\lambda)z_i\qquad\mbox{and}\qquad z_i\ge 0
	$$
	for $ i=1,\dots,k. $ Since $ f $ is matrix concave and non-negative we obtain
	$$
	f(\lambda y)\ge \lambda f(x)+(1-\lambda)f(z)\ge\lambda f(x)
	$$
	where $ z=(z_1,\dots,z_k). $	The result now follows by letting $ \lambda $ tend to one.
	\qed\vskip 2ex

	The converse is not true. The function of two variables $ f(r_1,r_2)=r_1r_2 $ is 
	indeed matrix increasing of any order in the sense that
	$$
	f(x_1,x_2)=x_1\otimes x_2\le y_1\otimes y_2=f(y_1,y_2)
	$$
	for $ 0\le x_1\le y_1 $ and $ 0\le x_2\le y_2, $ but it is not even concave as a 
	real function. However, the
	situation is quite different for functions of only one variable. 	Mathias~\cite{kn:mathias:1990} showed that a
	function, defined on the positive real half-line and matrix monotone of order $ n, $ is matrix
	concave of order $ [n/2]. $
	It follows from~\cite{kn:hansen:1980,kn:hansen:1982}, although not stated explicitely, that a function, defined 
	on the real positive half-line and matrix 
	monotone of order $ 4n, $ is matrix concave of order $ n. $ 
	If we relax Mathias' result only very slightly and are satisfied with proving that a
	function $ f:[0,\infty[\to{\mathbf R}, $ matrix monotone of order $ 2n, $ is matrix concave of 	order $ n, $
	then the following very simple argument will do. Let $ x_1,x_2 $ be positive definite matrices 	of order $ n $
	and notice~\cite{kn:hansen:1980} that to a given $ \varepsilon>0 $ the inequality
	$$
	V^*\left(\begin{array}{cc}
			    x_1 & 0\\
			    0   & x_2
			    \end{array}\right)V=\frac{1}{2}\left(\begin{array}{cc}
							                  x_1+x_2  & x_2-x_1\\
							                  x_2-x_1 & x_1+x_2
					                              \end{array}\right)
	\le\left(\begin{array}{cc}
		    2^{-1}(x_1+x_2)+\varepsilon & 0\\
		    0                           & \lambda
		    \end{array}\right)
	$$
	is valid for a sufficiently large $ \lambda>0, $  where
	$$
	V=\frac{\sqrt{2}}{2}\left(\begin{array}{cc}
	   			         1 & -1\\
					   1 & 1
					   \end{array}\right)
	$$
	is a unitary block matrix of order $ 2n\times 2n. $ We then obtain
	$$
	\begin{array}{l}\displaystyle
	\frac{1}{2}\left(\begin{array}{cc}
	                  f(x_1)+f(x_2) & f(x_2)-f(x_1)\\
	                  f(x_2)-f(x_1) & f(x_1)+f(x_2)
	                 \end{array}\right)=V^*\left(\begin{array}{cc}
			                                 f(x_1) & 0\\
								   0      & f(x_2)
							         \end{array}\right)V\\[4ex]
	=f\left(V^*\left(\begin{array}{cc}
				x_1 & 0\\
				0   & x_2
				\end{array}\right)V\right)\le\left(\begin{array}{cc}
		   								f(2^{-1}(x_1+x_2)+\varepsilon) & 0\\
								            0                              & f(\lambda)
		    								\end{array}\right)
	\end{array}	
	$$
	and consequently
	$$
	\frac{f(x_1)+f(x_2)}{2}\le f\left(\frac{x_1+x_2}{2}+\varepsilon\right)
	$$
	from which the statement follows by letting $ \varepsilon $ tend to zero, since matrix monotone functions
	of order greater or equal to two are continuous (even continuously differentiable).

	We shall finally mention that Korßnyi and others have considered a notion of operator monotonicity for functions
	of two variables defined on $ I^2 $ where $ I=]-1,1[. $ The notion is closely connected to the theory of 
	analytic functions of several variables, and in particular to a generalization of the Riesz-Herglotz 
	formula~\cite{kn:koranyi:1961,kn:singh:2001}. According to this theory the function
	$$
	g(r_1,r_2)=\frac{r_1r_2}{(1+r_1)(1+r_2)}\qquad r_1,r_2\in]0,1[
	$$
	would be called operator monotone, but this is not consistent with 
	Theorem~\ref{theorem:operator monotonicity and operator convexity} as the continuous function
	$$
	f(r_1,r_2)=\frac{r_1^2\,r_2^2}{(1+r_1)(1+r_2)}\qquad r_1,r_2\in[0,1[
	$$
	is not operator convex. Korßnyi's notion of operator monotononicity leads to no significant distinction between
	functions of one and two variables as does the theory of operator convex functions.

	\section{Decompositions and monotonicity}

	\begin{definition}\label{decomposition}\rm
	Let $ x $ be a positive invertible operator acting on a Hilbert space $ H. $
	We say that an $ l $-tuple $ (y_1,\dots,y_l) $ of positive invertible operators on $ H $ is a 
	decomposition of $ x $ (of length $ l) $ if 		
	\begin{equation}\label{equation:decomposition}
	y_1+\cdots+y_l=x.
	\end{equation}
	The $ l $-tuple $ a=(a_1,\dots,a_l) $ defined by setting $ a_i=x^{-1/2}y_i^{1/2} $ for $ i=1,\dots,l $
	is called the associated unitary row.
	\end{definition}		
	We recall~\cite{kn:araki:2000} that an $ l $-tuple $ a=(a_1,\dots,a_l) $ of operators 
	on a Hilbert space $ H $ is said to be a unitary row, if there exists a unitary operator 
	$ U $ on the direct
	sum of $ l $ copies of $ H $ such that $ (a_1,\dots,a_l) $ is the first row in the
	$ l\times l $ block matrix representation of $ U. $ The equation 
	\begin{equation}\label{sufficient condition for a unitary row}
	a_1a_1^*+\cdots+a_la_l^*=1\qquad\mbox{(the identity on $ H) $}
	\end{equation}
	is a necessary, but in general not sufficient condition for $ a=(a_1,\dots,a_l) $ to be a unitary row.

	The row $ a=(a_1,\dots,a_l) $ associated with the decomposition of $ x $ in the
	definition above satisfy condition~(\ref{sufficient condition for a unitary row}) since
	$$
	a_1a_1^*+\cdots+a_la_l^*=
	x^{-1/2}y_1x^{-1/2}+\cdots+x^{-1/2}y_lx^{-1/2}=1.
	$$
	Araki and the author proved that an $ l $-tuple $ a=(a_1,\dots,a_l) $ satisfying
	condition~(\ref{sufficient condition for a unitary row}) is a unitary row, if
	$ \dim\ker a_i=\dim\ker a_i^* $ for at least one $ i=1,\dots,l. $ The condition is
	trivially satisfied in this case since the operators $ a_1,\dots,a_l $ are invertible. The
	$ l $-tuple $ a=(a_1,\dots,a_l) $ in Definition~\ref{decomposition} is therefore indeed a unitary row.
	We notice that $ y_i=a_i^*xa_i $ for $ i=1,\dots,l. $

	\begin{definition}\label{definition of index}
	An index is a pair $ (l,j) $ of integers, where $ l\ge 2 $ and $ 0\le j\le l-1. $ 
	\end{definition}

	\begin{definition}\label{operator monotone function}\rm
	Let $ f:\,]0,\alpha_1[\times\cdots\times ]0,\alpha_k[\to{\mathbf R} $ be a 
	real function. The constants $ \alpha_1,\dots,\alpha_k $ may be plus infinity. 

	\begin{itemize}

	\item[$(i)$] We say that $ f $ is operator monotone of index $ (l,j), $ if $ f $ is
	continuous and
	$$
	\mbox{diag}\Bigl(f(y_{t_11},\dots,y_{t_kk})\Bigr)_{|t|=j\,(mod\,l)}
	\le f(x_1,\dots,x_k)L_{l^{k-1}}\leqno{(*)}
	$$
	for every $ k $-tuple $ x=(x_1,\dots,x_k) $ in the domain of $ f $ acting on any Hilbert
	spaces $ H_1,\dots,H_k $ and all decompositions
	$$
	y_{1i}+\cdots+y_{li}=x_i\qquad i=1,\dots,k
	$$
	where $ L_{l^{k-1}} $ is the $ l^{k-1}\times l^{k-1} $ block matrix with the unit
	operator on the tensor product $ H_1\otimes\cdots\otimes H_k $ in each entry. The index $ t $ is a multi-index of 
	the form $ t=(t_1,\dots,t_k), $ where $ t_i=1,\dots,l $ for $ i=1,\dots,k $ and weight 
	$ |t|=t_1+\cdots+t_k. $

	\item[$(ii)$] We say that $ f $ is matrix monotone of index $ (l,j) $ and order
	$ (n_1,\dots,n_k), $ if
	the same inequalities $(*)$ are satisfied for operators acting only on Hilbert
	spaces $ H_1,\dots,H_k $ of finite dimensions $ (n_1,\dots,n_k). $ 

	\end{itemize}
	\end{definition}

	It is not difficult to establish that a continuous function is operator monotone of index $ (l,j), $
	if and only if it is matrix monotone of index $ (l,j) $ and all orders $ (n_1,\dots,n_k). $ The proof follows a 	suggestion by L{\"o}wner as reported by Bendat and 
	Sherman~\cite[Lemma 2.2]{kn:bendat:1955}
	and can easily be adapted to the present situation.
	Furthermore, consider $ k $-tuples $ (m_1,\dots,m_k) $ and $ (n_1,\dots,n_k) $ such that $ m_i\le n_i $ for 
	$ i=1,\dots,k. $ If a function is matrix convex of order $ (n_1,\dots,n_k) $ then it is also matrix convex of order
	$ (m_1,\dots,m_k). $  Likewise, if a function
	is matrix monotone of index $ (l,j) $ and order $ (n_1,\dots,n_k), $ then it is also matrix monotone of index
	$ (l,j) $ and order $ (m_1,\dots,m_k). $

	\begin{proposition}
	A continuous real function $ f:\,]0,\alpha[\to{\mathbf R} $ is operator monotone of
	any given index $ (l,j), $ if and only if it is
	operator monotone. Likewise is $ f $ matrix monotone of any given index $ (l,j) $ 
	and order $ n, $ if and only if it is matrix monotone of order $ n. $
	\end{proposition} 

	\proof If we set $ k=1, $ the inequality $ (*) $ reads
	$$
	\begin{array}{llrl}
     f(y_{l1})\le f(x_1)\quad&\quad\mbox{for} &j& = 0\\[1ex]
	f(y_{j1})\le f(x_1)\quad&\quad\mbox{for} &j& = 1,\dots,l-1
	\end{array}
	$$
	where $ y_{11}+\cdots+y_{l1}=x_1 $ is a decomposition of $ x_1. $ These
	inequalities are trivially satisfied if $ f $ is operator monotone. If on the
	other hand one of the above inequalities are satisfied for a given index
	$ (l,j) $ and all decompositions
	of any $ x_1 $ in the domain of $ f, $ then $ f $ is operator monotone. The same
	reasoning applies to matrix monotone functions.
	\qed

	\begin{proposition}\label{an operator monotone function is separately operator monotone}
	Let $ f:\,]0,\alpha_1[\times\cdots\times ]0,\alpha_k[\to{\mathbf R} $ be a continuous function and consider
	for $ i=1,\dots,k $ the function of one variable
	$$
	g_i(r_i)=f(r_1,\dots,r_k)
	$$
	obtained from $ f $ by keeping all variables fixed except the ith variable. If $ f $ is operator monotone 
	of some index 
	$ (l,j), $ then $ g_i $ is operator monotone. Likewise, if $ f $ 	is matrix monotone of some index 
	$ (l,j) $ and order $ (n_1,\dots,n_k), $ then $ g_i $ is matrix monotone of order 
	$ n_i. $
	\end{proposition}

	\proof Let $ f $ be operator monotone (or matrix monotone) of some index $ (l,j) $ and assume $ i=1. $ 
	We choose operators $ y\le x $ in the domain of $ g_1. $ For some sufficiently small 
	$ \varepsilon>0 $ we set
	$$
	y_{m1}=\left\{\begin{array}{ll}
			   y &m=j+1\\[0.5ex]
			   \varepsilon+(x-y)/(l-1)\quad &m\ne j+1
		         \end{array}\right.\quad\mbox{and}\quad 
	y_{m2}=\left\{\begin{array}{ll}
			   r_2\quad &m=l-1\\[0.5ex]
			   \varepsilon\quad &m\ne l-1,
		          \end{array}\right.
	$$
	and for $ p=3,\dots,k $
	$$
	y_{mp}=\left\{\begin{array}{ll}
			   r_p\quad &m=l\\[0.5ex]	
			   \varepsilon\quad &m\ne l.
		          \end{array}\right.
	$$
	We thus have the decompositions $ y_{11}+\cdots+y_{l1}=x+(l-1)\varepsilon $ and
	$$ 
	y_{1p}+\cdots+y_{lp}=r_p+(l-1)\varepsilon\qquad p=2,\dots,l.
	$$
	By only considering the index $ t=(j+1,l-1,l,\dots,l) $
	with length $  |t|=j\,(mod\,l) $ in $(*),$ we obtain the inequality
	$$
	f(y,r_2,\dots,r_p)\le 
	f(x+(l-1)\varepsilon,r_2+(l-1)\varepsilon,\dots,r_l+(l-1)\varepsilon)
	$$
	from which the inequality $ g_1(y)\le g_1(x) $ is derived by letting $ \varepsilon $ tend to zero.
	\qed\vskip 2ex
	 
	To further investigate
	the content of Definition~\ref{operator monotone function} we set $ k=2 $ and $ l=2. $
	The inequality $ (*) $ then reads 
	$$
	\left(\begin{array}{cc}
		f(y_{11},y_{12}) & 0\\
		0                & f(y_{21},y_{22})
		\end{array}\right)
	\le\left(\begin{array}{cc}
		   f(x_1,x_2) & f(x_1,x_2)\\
		   f(x_1,x_2) & f(x_1,x_2)
		   \end{array}\right)\qquad j=0
	$$
	and
	$$
	\left(\begin{array}{cc}
		f(y_{11},y_{22}) & 0\\
		0                & f(y_{21},y_{12})
		\end{array}\right)
	\le\left(\begin{array}{cc}
		   f(x_1,x_2) & f(x_1,x_2)\\
		   f(x_1,x_2) & f(x_1,x_2)
		   \end{array}\right)\qquad j=1
	$$
	for decompositions $ y_{11}+y_{21}=x_1 $ and $ y_{12}+y_{22}=x_2. $ This
	is so because the solutions to the equation $ |t|=t_1+t_2=j\,(mod\,2) $ are the
	multi-indices $ (1,1),(2,2) $ for $ j=0 $ and $ (1,2),(2,1) $ for $ j=1. $ If we
	set $ k=2 $ and $ l=3 $ the inequality $ (*) $ reads
	$$
	\left(\begin{array}{ccc}
		f(y_{11},y_{22}) & 0                & 0\\
		0                & f(y_{21},y_{12}) & 0\\
	      0                & 0                & f(y_{31},y_{32}) 
		\end{array}\right)
	\le f(x_1,x_2)\left(\begin{array}{ccc}
		              1 & 1 & 1\\
		              1 & 1 & 1\\
	                    1 & 1 & 1
		              \end{array}\right)\qquad j=0
	$$
	for decompositions $ y_{11}+y_{21}+y_{31}=x_1 $ and 
	$ y_{12}+y_{22}+y_{32}=x_2. $ This is so because the solutions to the equation
	$ |t|=t_1+t_2=0\,(mod\,3) $ are the	multi-indices $ (1,2),(2,1),(3,3). $ Finally,
	if we set $ k=3 $ and $ l=2 $ the inequality $ (*) $ reads
	$$
	\begin{array}{rl}
	&\left(\begin{array}{cccc}
	f(y_{11},y_{12},y_{23}) & 0                       & 0
	& 0\\
	0                       & f(y_{11},y_{22},y_{13}) & 0
	& 0\\
      0                       & 0                       & f(y_{21},y_{12},y_{13}) 
	& 0\\
      0                       & 0                       & 0 
	& f(y_{21},y_{22},y_{23}) 
	\end{array}\right)\\[6ex]
	\le&f(x_1,x_2,x_3) L_4\qquad j=0
	\end{array}
	$$
	for decompositions $ y_{11}+y_{21}=x_1, $ $ y_{12}+y_{22}=x_2 $ and
	$ y_{13}+y_{23}=x_3. $ This is so because the solutions to the equation
	$ |t|=t_1+t_2+t_3=0\,(mod\,2) $ are the multi-indices $ (1,1,2),(1,2,1),(2,1,1) $
	and $ (2,2,2). $ 

	\begin{theorem}\label{theorem:ordinary convexity}
	Let $ f:\,]0,\alpha_1[\times\cdots\times]0,\alpha_k[\to{\mathbf R} $ be a 
	continuous, real function. If
	the function $ g:\,]0,\alpha_1[\times\cdots\times]0,\alpha_k[\to{\mathbf R} $
	defined by
	$$
	g(r_1,\dots,r_k)=r_1^{-1}\cdots r_k^{-1}\, f(r_1,\dots,r_k)
	$$ 
	is matrix monotone of some index $ (l,j) $ and order $ (l,\dots,l), $ then
	$ f $ is convex. 
	\end{theorem}
	
	\proof	
	We consider the simple root $ \beta=e^{2\pi i/l} $ of the polynomial 
	$ X^l-1 $ and set
	$$
	u=\mbox{diag}\left(\beta^p\right)_{p=1}^l
	$$
	which is a unitary matrix acting on $ {\mathbf C}^l. $ 
	We introduce projections
	$$
	P_j=(u^*)^j P u^j\qquad j=1,\dots,l
	$$
	where $ P $ defined by
	$$
	P=\frac{1}{l}\left(\begin{array}{ccc}
			            1\cdots 1\\
					  \vdots\ddots\vdots\\
					  1\cdots 1
					  \end{array}\right)
	$$
	is a one-dimensional projection acting on $ {\mathbf C^l}. $ We notice that
	$$
	P_j=\frac{1}{l}\left(\beta^{(q-p)j}\right)_{p,q=1}^l
	\qquad j=1,\dots,l
	$$
	and consequently
	$$
	\sum_{j=1}^l P_j=\frac{1}{l}\left(\sum_{j=1}^l \beta^{(q-p)j}\right)_{p,q=1}^l=E_l
	$$
	where $ E_l $ is the $ l\times l $ identity matrix. The projections $ P_1,\dots,P_l $ are 
	thus mutually orthogonal.
 	Let $ x_{1i},\dots,x_{li} $ be real numbers in $ ]0,\alpha_i[ $ and set
	$$
	x_i=\mbox{diag}\Bigl(x_{ji}\Bigr)_{j=1}^l\qquad i=1,\dots,k.
	$$  
	The $ k $-tuple $ ((1+l\varepsilon)x_1,\dots,(1+l\varepsilon)x_k) $ is for a
	sufficiently small $ \varepsilon>0 $ in the domain of $ g. $ We define
	$$
	y_{ji}=x_i^{1/2}(P_j+\varepsilon)x_i^{1/2}\qquad j=1,\dots l;\, i=1,\dots,k
	$$
	and calculate
	$$
	y_{1i}+\cdots+y_{li}=x_i^{1/2}\left(\sum_{j=1}^l (P_j+\varepsilon)\right) x_i^{1/2}
	=(1+l\varepsilon) x_i\qquad i=1,\dots,k.
	$$
	Since $ g $ is matrix monotone of index $ (l,j) $ and order $ (l,\dots,l) $ it
	follows that
	$$
	\mbox{diag}\Bigl(g(y_{t_1 1},\dots,y_{t_k k})\Bigr)_{|t|=j\,(mod\, l)}
	\le g((1+l\varepsilon)x_1,\dots,(1+l\varepsilon)x_k)L_{l^{k-1}}
	$$
	or inserting $ g(r_1,\dots,r_k)=r_1^{-1}\cdots r_k^{-1}f(r_1,\dots,r_k) $ that
	$$
	\begin{array}{l}
	\mbox{diag}\Bigl((y_{s_1 1}^{-1/2}\otimes\cdots\otimes y_{s_k k}^{-1/2})
	f(y_{s_1 1},\dots,y_{s_k k})(y_{s_1 1}^{-1/2}\otimes\cdots\otimes y_{s_k k}^{-1/2})
	\Bigr)_{|s|=j\,(mod\, l)}\\[2ex]
	\le c_\varepsilon^{-k}\hskip -1pt\left((x_1^{-1/2}\otimes\cdots\otimes x_k^{-1/2})
	f(c_\varepsilon x_1,\dots,c_\varepsilon x_k)
	(x_1^{-1/2}\otimes\cdots\otimes x_k^{-1/2})
	\right)_{|t|=|s|=j\,(mod\, l)}
	\end{array}
	$$
	where $ c_\varepsilon=1+l\varepsilon. $ Multiplying to the left and to the right with
	the self-adjoint matrix
	$$
	\mbox{diag}\Bigl(y_{s_1 1}^{1/2}\otimes\cdots\otimes y_{s_k k}^{1/2}
	\Bigr)_{|s|=j\,(mod\, l)}
	$$
	we obtain
	$$
	\begin{array}{l}
	\mbox{diag}\Bigl(f(y_{s_1 1},\dots,y_{s_k k})\Bigr)_{|s|=j\,(mod\, l)}\\[2ex]
	\le c_\varepsilon^{-k}
	\left((y_{t_1 1}^{1/2}x_1^{-1/2}\otimes\cdots\otimes y_{t_k k}^{1/2}x_k^{-1/2})
	f(c_\varepsilon x_1,\dots,c_\varepsilon x_k)\right.\times\\[1ex]
	\hskip 15em \left.(x_1^{-1/2}y_{s_1 1}^{1/2}\otimes\cdots\otimes x_k^{-1/2}y_{s_k k}^{1/2})
	\right)_{|t|=|s|=j\,(mod\, l)}.
	\end{array}
	$$
	We introduce for $ s_i=1,\dots,l $ and $ i=1,\dots,k $ the $ l\times l $ matrix
	$$
	Q_{s_i i}	=\frac{1}{x_{1i}+\cdots+x_{li}}\Bigl(x_{p i}^{1/2}x_{q i}^{1/2}\beta^{(q-p)s_i}
	\Bigr)_{p,q=1}^l.
	$$
	It is an easy calculation to show that $ Q_{s_i i} $ is a projection and that
	$$
	x_i^{1/2}P_{s_i}\,x_i^{1/2}=\frac{x_{1i}+\cdots+x_{li}}{l}Q_{s_i i}\qquad
	s_i=1,\dots,l;\,i=1,\dots,k.
	$$
	Multiplying the above inequality from the left and the right with the projection
	$$
	\mbox{diag}\Bigl(Q_{s_1 1}\otimes\cdots\otimes Q_{s_k k}\Bigr)_{|s|=j\,(mod\, l)}
	$$
	and letting $ \varepsilon $ tend to zero we thus obtain
	\begin{equation}\label{1th inequality}
	\begin{array}{l}\displaystyle           	f\left(\frac{x_{11}+\cdots+x_{l1}}{l},\dots,\frac{x_{1k}+\cdots+x_{lk}}{l}\right)
	\mbox{diag}\Bigl(Q_{s_1 1}\otimes\cdots\otimes Q_{s_k k}\Bigr)_{|s|=j\,(mod\, l)}\\[3ex]
	\hskip 0.5em\le\displaystyle
	\left(\frac{x_{11}+\cdots+x_{l1}}{l}\right)^{-1}\cdots
	\left(\frac{x_{1k}+\cdots+x_{lk}}{l}\right)^{-1}\times\\[2ex]
	\Bigl((x_1^{1/2}P_{t_1}\otimes\cdots\otimes x_k^{1/2}P_{t_k})	f(x_1,\dots,x_k)
	(P_{s_1}x_1^{1/2}\otimes\cdots\otimes P_{s_k}x_k^{1/2})\Bigr)_{|t|=|s|=j\,(mod\, l)}
	\end{array}
	\end{equation}
	where we used that 
	$$
	Q_{t_i i}\,y_{t_i i}^{1/2}x_i^{-1/2} \to\left(\frac{x_{1i}+\cdots+x_{li}}{l}\right)^{1/2}
 	Q_{t_i i}\,x_i^{-1/2} =\left(\frac{x_{1i}+\cdots+x_{li}}{l}\right)^{-1/2}x_i^{1/2}P_{t_i}
	$$
	as $ \varepsilon $ tends to zero. We notice that $(\ref{1th inequality})$ is an 
	$ l^{k-1}\times l^{k-1} $ block matrix inequality of $ l^k\times l^k $ matrices.
	Let us in order to examine the inequality calculate the entry
	$$
	\begin{array}{rl}
	&\left[x_1^{1/2} P_{t_1}\otimes\cdots\otimes x_k^{1/2} P_{t_k}\right]_{p q}=
	\left[x_1^{1/2} P_{t_1}\right]_{p_1 q_1}\cdots\left[x_k^{1/2} P_{t_k}\right]_{p_k q_k}\\[2ex]
	=&x_{p_1 1}^{1/2}\,	l^{-1}\beta^{(q_1-p_1)t_1}\cdots 
	x_{p_k k}^{1/2}\, l^{-1}\beta^{(q_k-p_k)t_k}=
	l^{-k} \beta^{(q-p)\cdot t}\, x_{p_1 1}^{1/2}\cdots x_{p_k k}^{1/2}
	\end{array}
	$$
	for $ p=(p_1,\dots,p_k) $ and $ q=(q_1,\dots,q_k) $ with 
	$ p_1,\dots,p_k,q_1,\dots,q_k=1,\dots,l. $ We proceed to calculate the entry
	$$
	\begin{array}{l}
	\left[(x_1^{1/2}P_{t_1}\otimes\cdots\otimes x_k^{1/2}P_{t_k})	f(x_1,\dots,x_k)
	(P_{s_1}x_1^{1/2}\otimes\cdots\otimes P_{s_k}x_k^{1/2})\right]_{pq}\\[2ex]
	=\displaystyle
	\sum_{u_1,\dots,u_k=1}^l l^{-k} \beta^{(u-p)\cdot t}\, x_{p_1 1}^{1/2}\cdots x_{p_k k}^{1/2}
 \left[f(x_1,\dots,x_k)(P_{s_1}x_1^{1/2}\otimes\cdots\otimes P_{s_k}x_k^{1/2})\right]_{uq}\\[2ex]
	=\displaystyle
	\sum_{u_1,\dots,u_k=1}^l l^{-k} \beta^{(u-p)\cdot t}\, x_{p_1 1}^{1/2}\cdots x_{p_k k}^{1/2}
f(x_{u_1 1},\dots,x_{u_k k})l^{-k}\beta^{(q-u)\cdot s} x_{q_1 1}^{1/2}\cdots x_{q_k k}^{1/2}\\[2ex]
	=l^{-2k}\beta^{q\cdot s-p\cdot t}
	x_{p_1 1}^{1/2}\cdots x_{p_k k}^{1/2}\,x_{q_1 1}^{1/2}\cdots x_{q_k k}^{1/2}
	\displaystyle\sum_{u_1,\dots,u_k=1}^l \beta^{(t-s)\cdot u} f(x_{u_1 1},\dots,x_{u_k k})	
	\end{array}
	$$
	where we used that $ f(x_1,\dots,x_k) $ is a diagonal matrix with 
	$ f(x_{u_1 1},\dots,x_{u_k k}) $ as the $u$th diagonal entry, and finally calculate the 	diagonal entry
	$$
	\left[Q_{q_1 1}\otimes\cdots\otimes Q_{q_k k}\right]_{qq}
	=\left[Q_{q_1 1}]_{q_1 q_1}\cdots\right[Q_{q_k k}]_{q_k q_k}
	=x_{q_1 1}\cdots x_{q_k k}\prod_{i=1}^k (x_{1i}+\cdots+x_{li})^{-1}.
	$$
	We obtain from $(\ref{1th inequality})$ an inequality between 
	$ l^{k-1}\times l^{k-1} $ matrices by
	retaining the $ (t,s) $-entry in each $ (t,s) $-block on both sides of the inequality
	and discarding all other entries. We then insert the entries calculated above in the 
	inequality obtained in this way and get
	$$
	\begin{array}{l}
	\displaystyle f\left(\frac{x_{11}+\cdots+x_{l1}}{l},\dots,\frac{x_{1k}+\cdots+x_{lk}}{l}\right)	\prod_{i=1}^k       	(x_{1i}+\cdots+x_{li})^{-1}\times\\[1ex]	
	\hfill\mbox{diag}\Bigl(x_{s_1 1}\cdots x_{s_k k}\Bigr)_{|s|=j\,(mod\, l)}\\[2ex]
	\le\displaystyle l^{-k}\prod_{i=1}^k (x_{1i}+\cdots+x_{li})^{-1}\times\\[1ex]
	\displaystyle
	\left(x_{t_1 1}^{1/2}\cdots x_{t_k k}^{1/2}x_{s_1 1}^{1/2}\cdots x_{s_k k}^{1/2}\hskip -0.4em
	\sum_{u_1,\dots,u_k=1}^l\hskip -0.3em\beta^{s\cdot s-t\cdot t+(t-s)\cdot u} f(x_{u_1 1},\dots,x_{u_k k})
	\right)_{|t|=|s|=j\,(mod\,l)}.
	\end{array}
	$$	
	Multiplying from the left and from the right with the self-adjoint matrix
	$$
	\mbox{diag}\left(x_{s_1 1}^{-1/2}\cdots x_{s_k k}^{-1/2}
	\prod_{i=1}^k (x_{1i}+\cdots+x_{li})^{1/2}\right)_{|s|=j\,(mod\,l)}	
	$$
	we obtain
	\begin{equation}\label{2th inequality}
	\begin{array}{l}\displaystyle
	f\left(\frac{x_{11}+\cdots+x_{l1}}{l},\dots,\frac{x_{1k}+\cdots+x_{lk}}{l}\right)
	E_{l^{k-1}}\\[2ex]
	\le\displaystyle l^{-k}\sum_{u_1,\dots,u_k=1}^l f(x_{u_1 1},\dots,x_{u_k k})
	\left(\beta^{s\cdot s-t\cdot t+(t-s)\cdot u}\right)_{|t|=|s|=j\,(mod\,l)}.
	\end{array}
	\end{equation}
	We define for each $ u=(u_1,\dots,u_k) $ with $ u_1,\dots,u_k=1,\dots,l $ an
	$ l^{k-1}\times l^{k-1} $ matrix $ \Pi_u $ by setting
	\begin{equation}\label{definition of projections}
	\Pi_u=l^{-(k-1)}\left(\beta^{s\cdot s-t\cdot t+(t-s)\cdot u}\right)_{|t|=|s|=j\,(mod\, l)}.
	\end{equation}
	It is an easy calculation to show that the matrices $ \Pi_u $ are self-adjoint projections,
	and the inequality $(\ref{2th inequality})$ can in terms of these projections be written as
	\begin{equation}\label{3th inequality}
	\begin{array}{l}\displaystyle
	f\left(\frac{x_{11}+\cdots+x_{l1}}{l},\dots,\frac{x_{1k}+\cdots+x_{lk}}{l}\right) E_{l^{k-1}}\\[2ex]
	\le\displaystyle\frac{1}{l}\sum_{u_1,\dots,u_k=1}^l f(x_{u_1 1},\dots,x_{u_k k})\Pi_u.
	\end{array}
	\end{equation}
	Because of
	$$
	\sum_{u_1,\dots,u_k=1}^l \beta^{(s-t)\cdot u}=l^k\delta_{ts}
	$$
	it follows that
	\begin{equation}\label{sum of projections}
	\sum_{u_1,\dots,u_k=1}^l \Pi_u=lE_{l^{k-1}}.
	\end{equation}
	Since each index $ (t,s) $ in each $ \Pi_u $ satisfy $ |t|=|s|=j\,(mod\,l), $ confer
	equation
	(\ref{definition of projections}), it follows that $ \Pi_u=\Pi_v $ for each $ v $ on the form 
	\begin{equation}\label{equivalent indices}
	v=(v_1,\dots,v_k)=(u_1+i\,(mod\,l),\dots,u_k+i\,(mod\,l))\qquad i=0,1,\dots,l-1.
	\end{equation}
	We also notice that for each $ u $ there are exactly $ l $ different indices in
	(\ref{equivalent indices}). 
	It follows that each projection is counted $ l $ times in the sum (\ref{sum of projections}).
	Two projections are consequently either orthogonal, or identical with their indices
	connected as in (\ref{equivalent indices}). Setting $ u=(1,\dots,1) $ and
	multiplying (\ref{3th inequality}) with $ \Pi_u $ we obtain
	$$
	\begin{array}{l}\displaystyle
	f\left(\frac{x_{11}+\cdots+x_{l1}}{l},\dots,\frac{x_{1k}+\cdots+x_{lk}}{l}\right)\Pi_u\\[2ex]
	\displaystyle\le\frac{1}{l}\Bigl(f(x_{11},\dots,x_{1k})+\cdots+f(x_{l1},\dots,x_{lk})\Bigr)\Pi_u.
	\end{array}
	$$
	Therefore $ f $ is convex.\qed\vskip 2ex
	
	A matrix monotone function may tend to minus infinity as the argument of the function
	approaches a point located on an axis, but it cannot go too fast.

	\begin{corollary}\label{corollary:growth of matrix monotone functions}
	Let $ g:\,]0,\alpha_1[\times\cdots\times]0,\alpha_k[\to{\mathbf R} $ be a 
	continuous function, which is matrix monotone of some index $ (l,j) $ and order
	$ (l,\dots,l). $ To each subset of the domain of $ g $ of the form 
	$ ]0,\beta_1[\times\cdots\times]0,\beta_k[ $ where $ \beta_1,\dots,\beta_k<\infty, $
	there is a constant $ C\ge 0 $ such that
	$$
	g(r_1,\dots,r_k)\ge-\frac{C}{r_1\cdots r_k}\qquad 
	(r_1,\dots,r_k)\in\,]0,\beta_1[\times\cdots\times]0,\beta_k[.
	$$	
	\end{corollary}

	\proof The function $ f:]0,\alpha_1[\times\cdots\times]0,\alpha_k[\to{\mathbf R} $
	given by
	$$
	f(r_1,\dots,r_k)=r_1\cdots r_k\, g(r_1,\dots,r_k)
	$$ 
	is convex by the preceeding theorem, and it is therefore bounded from below on bounded
	subsets of the domain.
	\qed\vskip 2ex

	To proceed, we need the following slight generalization 
	of~\cite[Theorem 1.2]{kn:araki:2000}. 
	\begin{theorem}\label{Jensen's operator inequality}
      Let $ f $ be a real, continuous function of $ k $ variables defined on the domain $ I_1\times\cdots\times I_k $
	where $ I_1,\dots,I_k $ are intervals containing zero and let $ (l,j) $ be any index. 
	The following statements are equivalent:

      \begin{itemize}

      \item[$(i)$] $ f $ is operator convex and
      $ f(r_1,\dots,r_k)\le 0 $ if $ r_i=0 $ for some $ i=1,\dots,k. $

      \item[$(ii)$] The operator inequality
	$$
      \begin{array}{rl}
      &\hskip -0.4em\mbox{diag}\Bigl(f(a_{s_11}^*x_1a_{s_11},\dots,a_{s_kk}^*x_ka_{s_kk})
      \Bigr)_{|s|=j\,(mod\,l)}\\[3ex]
      \le&\hskip -0.4em\Bigl((a_{t_11}^*\otimes\cdots\otimes a_{t_kk}^*)
      f(x_1,\dots,x_k)(a_{s_11}\otimes\cdots\otimes a_{s_kk})
      \Bigr)_{|t|=|s|=j\,(mod\,l)}
      \end{array}
      $$
	is valid for all unitary rows $ a_i=(a_{1i},\dots,a_{li}) $ of length $ l $
      acting on any Hilbert space $ H_i $ for $ i=1,\dots,k $ and all 
      $ k $-tuples $ (x_1,\dots,x_k) $ of
      self-adjoint operators in the domain of $ f $ acting on
      $ H_1,\dots,H_k. $
      
      \item[$(iii)$] The operator inequality
	$$
      \begin{array}{rl}
      &\hskip -0.4em\mbox{diag}\Bigl(f(p_{s_11}x_1p_{s_11},\dots,p_{s_kk}x_kp_{s_kk})
      \Bigr)_{|s|=j\,(mod\,l)}\\[3ex]
      \le&\hskip -0.4em\Bigl((p_{t_11}\otimes\cdots\otimes p_{t_kk})
      f(x_1,\dots,x_k)(p_{s_11}\otimes\cdots\otimes p_{s_kk})
      \Bigr)_{|t|=|s|=j\,(mod\,l)}
      \end{array}
      $$
	is valid for all partitions of unity
      $ p_{1i}+\cdots+p_{li}=1 $ on any Hilbert space $ H_i $ by orthogonal projections for each
      $ i=1,\dots,k $ and
      all $ k $-tuples $ (x_1,\dots,x_k) $ of
      self-adjoint operators in the domain of $ f $ acting on
      $ H_1,\dots,H_k. $
     
      \end{itemize}
	The indices $ s,t $ in $ (ii) $ and $ (iii) $ are multi-indices of the form 
      $ s=(s_1,\dots,s_k), $ where $ s_i=1,\dots,l $ for $ i=1,\dots,k $ 
      with weight $ |s|=s_1+\cdots+s_k. $

      \end{theorem}

  	In the reference~\cite{kn:araki:2000} the sufficiency of $ (ii) $ and $ (iii) $
	in order to obtain $ (i) $ were only established for indices of the form $ (l,0). $  
	However, rewriting of the original proof shows, mutatis mutandis, that the inequalities
	are indeed sufficient for the operator convexity of $ f $ for any index. The theorem above is
	stated for more general domains of the
	function $ f $ than in the original reference, cf. the discussion in the survey article~\cite{kn:hansen:2000:2}.
	It has the following version for functions of matrices~\cite{kn:hansen:2000:2}.

	\begin{theorem}\label{Jensen's matrix inequality}
      Let $ f $ be a real, continuous function of $ k $ variables defined on the domain $ I_1\times\cdots\times I_k $
	where $ I_1,\dots,I_k $ are intervals containing zero and let $ (l,j) $ be any index. 
	Let $ (n_1,\dots,n_k) $ be a $ k $-tuple of natural numbers and consider the statements:
 
      \begin{itemize}

      \item[$(i)$] $ f $ is matrix convex of order $ (ln_1,\dots,ln_k) $ and
      $ f(r_1,\dots,r_k)\le 0 $ if $ r_i=0 $ for some $ i=1,\dots,k. $

      \item[$(ii)$] The matrix inequality
	$$
      \begin{array}{rl}
      &\hskip -0.4em\mbox{diag}\Bigl(f(a_{s_11}^*x_1a_{s_11},\dots,a_{s_kk}^*x_ka_{s_kk})
      \Bigr)_{|s|=j\,(mod\,l)}\\[3ex]
      \le&\hskip -0.4em\Bigl((a_{t_11}^*\otimes\cdots\otimes a_{t_kk}^*)
      f(x_1,\dots,x_k)(a_{s_11}\otimes\cdots\otimes a_{s_kk})
      \Bigr)_{|t|=|s|=j\,(mod\,l)}
      \end{array}
      $$
	is valid for all unitary rows $ a_i=(a_{1i},\dots,a_{li}) $ of length $ l $
      acting on a Hilbert space $ H_i $ of dimension $ n_i $ for $ i=1,\dots,k $ and all 
      $ k $-tuples $ (x_1,\dots,x_k) $ of self-adjoint operators in the domain of $ f $ acting on
      $ H_1,\dots,H_k. $
      
      \item[$(iii)$] The matrix inequality
	$$
      \begin{array}{rl}
      &\hskip -0.4em\mbox{diag}\Bigl(f(p_{s_11}x_1p_{s_11},\dots,p_{s_kk}x_kp_{s_kk})
      \Bigr)_{|s|=j\,(mod\,l)}\\[3ex]
      \le&\hskip -0.4em\Bigl((p_{t_11}\otimes\cdots\otimes p_{t_kk})
      f(x_1,\dots,x_k)(p_{s_11}\otimes\cdots\otimes p_{s_kk})
      \Bigr)_{|t|=|s|=j\,(mod\,l)}
      \end{array}
      $$
	is valid for all partitions of unity
      $ p_{1i}+\cdots+p_{li}=1 $ on a Hilbert space $ H_i $ of dimension $ n_i $ by orthogonal projections for each
      $ i=1,\dots,k $ and all $ k $-tuples $ (x_1,\dots,x_k) $ of
      self-adjoint operators in the domain of $ f $ acting on
      $ H_1,\dots,H_k. $

	\item[$(iv)$] The matrix inequality
	$$
      \begin{array}{rl}
      &\hskip -0.4em\mbox{diag}\Bigl(f(p_{s_11}x_1p_{s_11},\dots,p_{s_kk}x_kp_{s_kk})
      \Bigr)_{|s|=j\,(mod\,l)}\\[3ex]
      \le&\hskip -0.4em\Bigl((p_{t_11}\otimes\cdots\otimes p_{t_kk})
      f(x_1,\dots,x_k)(p_{s_11}\otimes\cdots\otimes p_{s_kk})
      \Bigr)_{|t|=|s|=j\,(mod\,l)}
      \end{array}
      $$
	is valid for all partitions of unity
      $ p_{1i}+\cdots+p_{li}=1 $ on a Hilbert space $ H_i $ of dimension $ ln_i $ by orthogonal projections for each
      $ i=1,\dots,k $ and all $ k $-tuples $ (x_1,\dots,x_k) $ of
      self-adjoint operators in the domain of $ f $ acting on
      $ H_1,\dots,H_k. $

	\item[$(v)$] $ f $ is matrix convex of order $ (n_1,\dots,n_k) $ and
      $ f(r_1,\dots,r_k)\le 0 $ if $ r_i=0 $ for some $ i=1,\dots,k. $
     
      \end{itemize}

	\noindent The implications $ (i)\Rightarrow (ii)\Rightarrow (iii) $ and $ (iv)\Rightarrow (v) $ are then valid. 

      \end{theorem}

	The indices $ s,t $ in $ (ii),(iii) $ and $ (iv) $ are multi-indices of the form 
      $ s=(s_1,\dots,s_k), $ where $ s_i=1,\dots,l $ for $ i=1,\dots,k $ 
      with weight $ |s|=s_1+\cdots+s_k. $ Since the Hilbert spaces in $ (ii) $ are finite dimensional,
	it follows that any row $ a_i=(a_{1i},\dots,a_{li}) $ satisfying 
	condition~(\ref{sufficient condition for a unitary row}) is unitary.

	\begin{theorem}\label{convexity of f implies monotonicity of g}
	Let $ f:[0,\alpha_1[\times\cdots\times [0,\alpha_k[\to{\mathbf R} $ be a continuous real function
	such that $ f(r_1,\dots,r_k)\le 0 $ if $ r_i=0 $ for some $ i=1,\dots,k. $ The constants 
	$ \alpha_1,\dots,\alpha_k $ may be plus infinity. If $ f $ is matrix convex of order $ (ln_1,\dots,ln_k) $
	for some integer $ l\ge 2 $ and some $ k $-tuple of natural numbers $ (n_1,\dots,n_k), $ 
	then the function
	$$
	g(r_1,\dots,r_k)=r_1^{-1}\cdots r_k^{-1} f(r_1,\dots,r_k)
	\qquad (r_1,\dots,r_k)\in ]0,\alpha_1[\times\cdots\times ]0,\alpha_k[ 
	$$
	is matrix monotone of index $ (l,j) $ and order $ (n_1,\dots,n_k) $
	for $ j=0,1,\dots,l-1. $
	\end{theorem}

	\proof  Let $ (x_1,\dots,x_k) $ be any $ k $-tuple of positive invertible operators 
	in the domain of $ f $ acting on Hilbert spaces
	$ H_1,\dots,H_k $ of dimensions $ n_1,\dots,n_k $ and let 
	$$
	y_{1i}+\cdots+y_{li}=x_i
	$$
	be any decomposition of $ x_i $ of length $ l $ for each $ i=1,\dots,k. $ We set
	$$
	a_{s_ii}=x_i^{-1/2}y_{s_ii}^{1/2}\qquad	s_i=1,\dots,l;\, i=1,\dots,k
	$$
	and observe that
	$$
	y_{s_ii}=a_{s_ii}^*x_ia_{s_ii}\qquad s_i=1,\dots,l;\, i=1,\dots,k.
	$$
	If $ f $ is matrix convex of order $ (ln_1,\dots,ln_k) $ we may apply
	Jensen's matrix inequality for functions of several variables, cf. Theorem~\ref{Jensen's matrix inequality}
	$ (i)\Rightarrow (ii). $ For each $ j=0,1,\dots,l-1 $ we have
	$$
	\begin{array}{rl}
	&\mbox{diag}\Bigl(f(y_{s_11},\dots,y_{s_kk})\Bigr)_{|s|=j\,(mod\,l)}\\[3ex]
	=&\mbox{diag}\Bigl(f(a_{s_11}^*x_1a_{s_11},\dots,a_{s_kk}^*x_ka_{s_kk})
      \Bigr)_{|s|=j\,(mod\,l)}\\[3ex]
      \le&\hskip -0.4em\Bigl((a_{t_11}^*\otimes\cdots\otimes a_{t_kk}^*)
      f(x_1,\dots,x_k)(a_{s_11}\otimes\cdots\otimes a_{s_kk})
      \Bigr)_{|t|=|s|=j\,(mod\,l)}\\[3ex]
	=&\hskip -0.4em\Bigl((y_{t_11}^{1/2}x_1^{-1/2}\otimes\cdots\otimes y_{t_kk}^{1/2}x_k^{-1/2})
      f(x_1,\dots,x_k)\times\\[1ex]
	&\hfill (x_1^{-1/2}y_{s_11}^{1/2}\otimes\cdots\otimes x_k^{-1/2}y_{s_kk}^{1/2})
      \Bigr)_{|t|=|s|=j\,(mod\,l)}\\[3ex]
	=&\hskip -0.4em\Bigl((y_{t_11}^{1/2}\otimes\cdots\otimes y_{t_kk}^{1/2})
      g(x_1,\dots,x_k)(y_{s_11}^{1/2}\otimes\cdots\otimes y_{s_kk}^{1/2})
      \Bigr)_{|t|=|s|=j\,(mod\,l)}
      \end{array}
      $$
	and multiplying to the left and to the right with the self-adjoint operator
	$$
	C=\mbox{diag}\Bigl(y_{s_11}^{-1/2}\otimes\cdots\otimes y_{s_kk}^{-1/2}\Bigr)_{|s|=j\,(mod\,l)}
	$$
	in the above inequality, we obtain
	$$
	\mbox{diag}\Bigl((y_{s_11}^{-1}\otimes\cdots\otimes
	y_{s_kk}^{-1})f(y_{s_11},\dots,y_{s_kk})\Bigr)_{|s|=j\,(mod\,l)}\le
	\Bigl(g(x_1,\dots,x_k)\Bigr)_{|t|=|s|=j\,(mod\,l)}
	$$
	or equivalently
	$$
	\mbox{diag}\Bigl(g(y_{s_11},\dots,y_{s_kk})\Bigr)_{|s|=j\,(mod\,l)}\le g(x_1,\dots,x_k)L_{l^{k-1}}
	$$
	showing that $ g $ is matrix monotone of index $ (l,j) $ and order $ (n_1,\dots,n_k). $
	\qed

	\begin{theorem}\label{monotonicity of g implies convexity of f}
	Let $ f:[0,\alpha_1[\times\cdots\times [0,\alpha_k[\to{\mathbf R} $ be a continuous, 
	real function and suppose the function
	$$
	g(r_1,\dots,r_k)=r_1^{-1}\cdots r_k^{-1}f(r_1,\dots,r_k)
	\qquad (r_1,\dots,r_k)\in\, ]0,\alpha_1[\times\cdots\times ]0,\alpha_k[ 
	$$
	is matrix monotone of some index $ (l,j) $ and order $ (ln_1,\dots,ln_k). $ Then the 
	following statements are valid:

	\begin{itemize}
	
	\item[$(i)$] $ f(r_1,\dots,r_k)\le 0 $ if $ r_i=0 $ for some $ i=1,\dots,k. $

	\item[$(ii)$] $ f $ is matrix convex 
	of order $ (n_1,\dots,n_k). $

	\end{itemize}
	\end{theorem}

	\proof Since $ g $ is an increasing function in each coordinate, 
	cf. Proposition \ref{an operator monotone function is separately operator monotone}, 
	the first statement follows.

	Let $ (x_1,\dots,x_k) $ be a $ k $-tuple of positive invertible operators in the domain 	of $ f $
	acting on Hilbert spaces $ H_1,\dots,H_k $ of dimensions $ ln_1,\dots,ln_k $ and let
	$$
	p_{1i}+\cdots+p_{li}=1\qquad i=1,\dots,k
	$$ 
	be resolutions of the identity on $ H_i $ of length $ l. $ We choose a positive 
	$ \varepsilon $ such that $ (1+l\varepsilon)x $ is in the domain of $ f $ and set
	$$
	y_{s_i i}=x_i^{1/2}(p_{s_i i}+\varepsilon)x_i^{1/2}\qquad s_i=1,\dots,l;\, i=1,\dots,k.
	$$
	We consider the decompositions
	$$
	y_{1i}+\cdots+y_{li}=(1+l\varepsilon)x_i\qquad i=1,\dots,k
	$$
	and use the assumption to obtain
	$$
	\mbox{diag}\Bigl(g(y_{s_1 1},\dots,y_{s_k k})\Bigr)_{|s|=j\,(mod\,l)}
	\le g((1+l\varepsilon)x_1,\dots,(1+l\varepsilon)x_k)L_{l^{k-1}}.
	$$
	We introduce the diagonal block matrix
	$$
	C=\mbox{diag}\Bigl(x_1^{1/2}(p_{s_1 1}+\varepsilon)x_1 p_{s_1 1}\otimes\cdots\otimes
	x_k^{1/2}(p_{s_k k}+\varepsilon)x_k p_{s_k k}\Bigr)_{|s|=j\,(mod\,l)}
	$$
	and multiply to the left with $ C^* $ and to the right with $ C $ in the above inequality to obtain
	$$
	\begin{array}{l}
	\begin{array}{l}
	\mbox{diag}\Bigl(
	(p_{s_1 1}x_1(p_{s_1 1}+\varepsilon)x_1^{1/2}\otimes\cdots\otimes p_{s_k k}x_k(p_{s_k k}+\varepsilon)x_k^{1/2})
	g(y_{s_1 1},\dots,y_{s_k k})\times\\[1ex]
	\hfill
	(x_1^{1/2}(p_{s_1 1}+\varepsilon)x_1 p_{s_1 1}\otimes\cdots\otimes x_k^{1/2}(p_{s_k k}+\varepsilon)x_k p_{s_k k})
	\Bigr)_{|s|=j\,(mod\,l)}
	\end{array}\\[3ex]
	\begin{array}{l}
	\le\Bigl(
	(p_{t_1 1}x_1(p_{t_1 1}+\varepsilon)x_1^{1/2}\otimes\cdots\otimes p_{t_k k}x_k(p_{t_k k}+\varepsilon)x_k^{1/2}) 	g((1+l\varepsilon)x_1,\dots,(1+l\varepsilon)x_k)\\[1ex]
	\hfil\times
	(x_1^{1/2}(p_{s_1 1}+\varepsilon)x_1 p_{s_1 1}\otimes\cdots\otimes x_k^{1/2}(p_{s_k k}+\varepsilon)x_k p_{s_k k})
	\Bigr)_{|t|=|s|=j\,(mod\,l)}.
	\end{array}
	\end{array}
	$$
	Inserting 
	$$
	\begin{array}{l}
	g(y_{s_1 1},\dots,y_{s_k k})=
	(x_1^{-1/2}(p_{s_1 1}+\varepsilon)^{-1}x_1^{-1/2}\otimes\cdots\otimes
	x_k^{-1/2}(p_{s_k k}+\varepsilon)^{-1}x_k^{-1/2})\times\\[1ex]
	\hfill f(x_1^{1/2}(p_{s_1 1}+\varepsilon)x_1^{1/2},\dots,x_k^{1/2}(p_{s_k k}+\varepsilon)x_k^{1/2})
	\end{array}
	$$
	and 
	$$
	\begin{array}{l}
	g((1+l\varepsilon)x_1,\dots,(1+l\varepsilon)x_k)=\\[1ex]
	(1+l\varepsilon)^{-k}
	(x_1^{-1/2}\otimes\cdots\otimes x_k^{-1/2})f((1+l\varepsilon)x_1,\dots,(1+l\varepsilon)x_k)
	(x_1^{-1/2}\otimes\cdots\otimes x_k^{-1/2})
	\end{array}
	$$
	in the inequality, and then letting $ \varepsilon $ tend to zero we obtain
	$$
	\begin{array}{l}
	\begin{array}{l}
	\mbox{diag}\Bigl(
	(p_{s_1 1}x_1^{1/2}\otimes\cdots\otimes p_{s_k k}x_k^{1/2})
	f(x_1^{1/2}p_{s_1 1}x_1^{1/2},\dots,x_k^{1/2}p_{s_k k}x_k^{1/2})\times\\[1ex]
	\hfill
	(x_1^{1/2}p_{s_1 1}x_1 p_{s_1 1}\otimes\cdots\otimes x_k^{1/2}p_{s_k k}x_k p_{s_k k})
	\Bigr)_{|s|=j\,(mod\,l)}
	\end{array}\\[3ex]
	\le\Bigl((p_{t_1 1}x_1 p_{t_1 1}\otimes\cdots\otimes p_{t_k k}x_k p_{t_k k})f(x_1,\dots,x_k)\times\\[1ex]
	\hfill(p_{s_1 1}x_1 p_{s_1 1}\otimes\cdots\otimes p_{s_k k}x_k p_{s_k k})
	\Bigr)_{|t|=|s|=j\,(mod\,l)}.
	\end{array}
	$$ 
	The identity
	$$
	\begin{array}{l}
	f(x_1^{1/2}p_{s_1 1}x_1^{1/2},\dots,x_k^{1/2}p_{s_k k}x_k^{1/2})
	(x_1^{1/2}p_{s_1 1}\otimes\cdots\otimes x_k^{1/2}p_{s_k k})=\\[1ex]
	\hskip 10em (x_1^{1/2} p_{s_1 1}\otimes\cdots\otimes x_k^{1/2}p_{s_k k})
	f(p_{s_1 1}x_1 p_{s_1 1},\dots,p_{s_k k}x_k p_{s_k k})
	\end{array}
	$$
	follows by first considering polynomials and then applying Weierstrass' approximation theorem. Inserting
	the identity in the inequality above we obtain
	$$
	\begin{array}{l}
	\begin{array}{l}
	\mbox{diag}\Bigl(
	(p_{s_1 1}x_1 p_{s_1 1}\otimes\cdots\otimes p_{s_k k}x_k p_{s_k k})
	f(p_{s_1 1}x_1 p_{s_1 1},\dots,p_{s_k k}x_k p_{s_k k})\times\\[1ex]
	\hfill
	(p_{s_1 1}x_1 p_{s_1 1}\otimes\cdots\otimes p_{s_k k}x_k p_{s_k k})
	\Bigr)_{|s|=j\,(mod\,l)}
	\end{array}\\[3ex]
	\le\Bigl((p_{t_1 1}x_1 p_{t_1 1}\otimes\cdots\otimes p_{t_k k}x_k p_{t_k k})f(x_1,\dots,x_k)\times\\[1ex]
	\hfill (p_{s_1 1}x_1 p_{s_1 1}\otimes\cdots\otimes p_{s_k k}x_k p_{s_k k})
	\Bigr)_{|t|=|s|=j\,(mod\,l)}
	\end{array}
	$$ 
	and hence
	$$
	\begin{array}{l}
	\mbox{diag}\Bigl(
	(p_{s_1 1}\otimes\cdots\otimes p_{s_k k}) f(p_{s_1 1}x_1 p_{s_1 1},\dots,p_{s_k k}x_k p_{s_k k})
	(p_{s_1 1}\otimes\cdots\otimes p_{s_k k})
	\Bigr)_{|s|=j\,(mod\,l)}\\[2ex]	
	\le\Bigl((p_{t_1 1}\otimes\cdots\otimes p_{t_k k})f(x_1,\dots,x_k)(p_{s_1 1}\otimes\cdots\otimes p_{s_k k})
	\Bigr)_{|t|=|s|=j\,(mod\,l)}.
	\end{array}
	$$ 
	Because of $ (i) $ we obtain
	$$
	\begin{array}{l}
	f(p_{s_1 1}x_1 p_{s_1 1},\dots,p_{s_k k}x_k p_{s_k k})\\[2ex]
	\le
	(p_{s_1 1}\otimes\cdots\otimes p_{s_k k}) f(p_{s_1 1}x_1 p_{s_1 1},\dots,p_{s_k k}x_k p_{s_k k})
	(p_{s_1 1}\otimes\cdots\otimes p_{s_k k})
	\end{array}
	$$
	and consequently
	$$
	\begin{array}{l}
	\mbox{diag}\Bigl(f(p_{s_1 1}x_1 p_{s_1 1},\dots,p_{s_k k}x_k p_{s_k k})\Bigr)_{|s|=j\,(mod\,l)}\\[2ex]	
	\le\Bigl((p_{t_1 1}\otimes\cdots\otimes p_{t_k k})f(x_1,\dots,x_k)(p_{s_1 1}\otimes\cdots\otimes p_{s_k k})
	\Bigr)_{|t|=|s|=j\,(mod\,l)}
	\end{array}
	$$ 
	which is Jensen's matrix inequality. We thus deduce, cf. Theorem~\ref{Jensen's matrix inequality}
	$ (iv)\Rightarrow (v), $ that $ f $ is matrix convex of order
	$ (n_1,\dots,n_k). $\qed\vskip 2ex

	One may think that the preceeding theorem, which ensures matrix convexity of 
	$ f, $ could replace Theorem~\ref{theorem:ordinary convexity} which with
	similar conditions only imparts ordinary convexity on $ f. $ However, it
	is essential in the proof of the preceeding theorem that $ f $ is defined also
	on the axes, while this is not required in Theorem~\ref{theorem:ordinary convexity}.
	This problem can easily be overcome for functions of only one variable by making a small
	translation of the matrix monotone function $ g. $ This remedy is not available for functions 
	of several variables, since the translation of the decomposition of an operator no
	longer is a decomposition of the translated operator, cf. equation~(\ref{equation:decomposition}).

	\begin{corollary}
	Let $ g:\,]0,\alpha_1[\times\cdots\times ]0,\alpha_k[\to{\mathbf R} $ be a 
	continuous real function.
	If $ g $ is matrix monotone of some index $ (l,j) $ and order $ (lm n_1,\dots,lm n_k) $ 
	for a natural number $ m $ and a $ k $-tuple of natural numbers $ (n_1,\dots,n_k), $
	then it is matrix monotone of index
 	$ (m,h) $ and order $ (n_1,\dots,n_k) $ for $ h=0,1,\dots,m-1. $
	\end{corollary}

	\proof The real and continuous function $ f $ defined by
	$$
	f(r_1,\dots,r_k)=r_1\cdots r_k\, g(r_1,\dots,r_k)\qquad 
	0 < r_i < \alpha_i\quad\mbox{for}\quad i=1,\dots,k
	$$
	is convex by Theorem~\ref{theorem:ordinary convexity}. It therefore extends
	to a continuous function
	$$
	\tilde f:[0,\alpha_1[\times\cdots\times[0,\alpha_k[\to{\mathbf R},
	$$
	and it follows that $ \tilde f(r_1,\dots,r_k)\le 0 $ if $ r_i=0 $ for some $ i=1,\dots,k. $
	The function $ \tilde f $ is matrix convex of order
	$ (mn_1,\dots,mn_k) $ by Theorem~\ref{monotonicity of g implies convexity of f}. 
	The function $ g $
	is thus matrix monotone of index $ (m,h) $ and order $ (n_1,\dots,n_k) $
	for $ h=0,1,\dots,m-1 $ by 
	Theorem \ref{convexity of f implies monotonicity of g}.\qed

	\begin{corollary}\label{independency of a particular index}
	Let $ g:\,]0,\alpha_1[\times\cdots\times ]0,\alpha_k[\to{\mathbf R} $ be a 
	continuous function.
	If $ g $ is operator monotone of some index, then it is operator monotone of all indices.
	\end{corollary}

	\begin{definition}\label{ex post definition of operator monotonicity} We say that a continuous function 
	$ g:\,]0,\alpha_1[\times\cdots\times ]0,\alpha_k[\to{\mathbf R} $ is operator monotone, if it is operator
	monotone of some and hence operator monotone of all indices.
	\end{definition}

	\noindent{\bf Proof} (of Theorem~\ref{theorem:operator monotonicity and operator convexity}):\hskip 0.5em
	The statement follows
	by combining Theorem~\ref{convexity of f implies monotonicity of g}, 
	Theorem~\ref{monotonicity of g implies convexity of f} and 
	Definition~\ref{ex post definition of operator monotonicity}.
	\qed\vskip 2ex
	
	The simplest example of operator convex functions satisfying the boundary conditions
	in Theorem~\ref{theorem:operator monotonicity and operator convexity} are the negative
	constants. The function
	$$
	g(r_1,\dots,r_k)=-r_1^{-1}\cdots r_k^{-1}
	$$
	defined in the first (open) quadrant is thus operator monotone, cf. also 
	Corollary~\ref{corollary:growth of matrix monotone functions}. The set of operator monotone
	functions defined on a given domain is a weakly closed convex cone, but
	the constant function $ g(r_1,\dots,r_k)=1 $ is not operator monotone for
	$ k\ge 2. $ This must indeed be so since the function $ (r_1,r_2)\to r_1r_2 $ is not convex.

    	\bibliographystyle{c:/pctexv4/texinput/plain}
	\bibliography{c:/unidata/macros/mathharv}	

\begin{thebibliography}{10}

\bibitem{kn:araki:2000}
H.~Araki and F.~Hansen.
\newblock Jensen's operator inequality for functions of several variables.
\newblock {\em Proc. Amer. Math. Soc.}, 128:2075--2084, 2000.

\bibitem{kn:bendat:1955}
J.~Bendat and S.~Sherman.
\newblock Monotone and convex operator functions.
\newblock {\em Trans. Amer. Math. Soc.}, 79:58--71, 1955.

\bibitem{kn:hansen:1980}
F.~Hansen.
\newblock An operator inequality.
\newblock {\em Math. Ann.}, 246:249--250, 1980.

\bibitem{kn:hansen:1997:2}
F.~Hansen.
\newblock Operator convex functions of several variables.
\newblock {\em Publ. RIMS, Kyoto Univ.}, 33:443--463, 1997.

\bibitem{kn:hansen:2000:2}
F.~Hansen.
\newblock Operator inequalities associated with \uppercase{J}ensen's
  inequality.
\newblock In T.M. Rassias, editor, {\em Survey on Classical Inequalities},
  pages 67--98. Kluwer Academic Publishers, 2000.

\bibitem{kn:hansen:2002}
F.~Hansen.
\newblock Convex trace functions of several variables.
\newblock {\em Linear Algebra Appl.}, 341:309--315, 2002.

\bibitem{kn:hansen:1982}
F.~Hansen and G.K. Pedersen.
\newblock Jensen's inequality for operators and \uppercase{L}{\"o}wner's
  theorem.
\newblock {\em Math. Ann.}, 258:229--241, 1982.

\bibitem{kn:kadison:1982/1986}
R.V. Kadison and J.R. Ringrose.
\newblock {\em Fundamentals of the Theory of Operator Algebras, I-II}.
\newblock {Academic Press}, 1982 \& 1986.

\bibitem{kn:koranyi:1961}
A.~Kor{\'a}nyi.
\newblock On some classes of analytic functions of several variables.
\newblock {\em Trans Amer. Math. Soc.}, 101:520--554, 1961.

\bibitem{kn:Lieb:2001}
E.~Lieb and G.K. Pedersen.
\newblock Convex multivariate trace functions.
\newblock {\em Rev. Math. Phys., to appear}, 2002.

\bibitem{kn:mathias:1990}
R.~Mathias.
\newblock Concavity of monotone matrix functions of finite order.
\newblock {\em Linear and Multilinear Algebra}, 27:129--138, 1990.

\bibitem{kn:pedersen:2001}
G.K. Pedersen.
\newblock Convex trace functions of several variables.
\newblock {\em Preprint}, 2001.

\bibitem{kn:singh:2001}
M.~Singh and H.L. Vasudeva.
\newblock Monotone matrix functions of two variables.
\newblock {\em Linear Algebra Appl.}, 328:131--152, 2001.

\end{thebibliography}
	          
      \nocite{kn:hansen:1980}
      \nocite{kn:hansen:1982}
	\nocite{kn:hansen:2002}
      \nocite{kn:kadison:1982/1986}
          
      \vfill

      {\small\noindent Frank Hansen: Institute of Economics, University
       of Copenhagen, Studiestraede 6, DK-1455 Copenhagen K, Denmark.}

      \end{document}